\documentstyle[11pt,epsfig]{amsart}
\theoremstyle{plain}
\def\CM{\mathcal M}

\def\IZ{\mathbb Z}

\def\a{\alpha}
\def\b{\beta}

\def\S{\Sigma}

\title[]{Minimal Number of Singular Fibers in a Lefschetz Fibration }
\author{Mustafa Korkmaz}
\address{Department of Mathematics \\
Middle East Technical University \\ 06531 Ankara, T\"urk\.iye}
\author{Burak Ozbagci}
\address{Department of Mathematics\\
University of California Irvine, CA 92697}

\begin{document}
\newtheorem{thm}{Theorem}
\newtheorem{cor}[thm]{Corollary}
\newtheorem{prop}[thm]{Proposition}
\newtheorem{lemma}[thm]{Lemma}
\newtheorem{df}[thm]{Definition}
\newtheorem{thmo}{Theorem}
\maketitle
\setcounter{secnumdepth}{2} \setcounter{section}{0}

\maketitle \setcounter{section}{0}


\begin{quote}
 ABSTRACT: There exists a (relatively minimal) genus $g$ Lefschetz
fibration with only one singular fiber over a closed (Riemann)
surface of genus $h$ iff $g \geq 3$ and $h \geq 2$. The singular
fiber can be chosen to be reducible or irreducible. Other results
are that every Dehn twist on a closed surface of genus at least
three is a product of two commutators and no Dehn twist on any
closed surface is equal to a single commutator.
\end{quote}

\section{Introduction}

A Lefschetz fibration is a fibration of a smooth 4-manifold over a 
surface with general fiber another closed orientable surface, which 
may admit certain singular fibers. The isomorphism class of a Lefschetz 
fibration is determined by its global monodromy. This leads to
a combinatorial approach to study the topology of smooth 4-manifolds 
(which admit Lefschetz fibrations) by their monodromy representations in the
mapping class groups. In particular, the number of singular fibers 
in a Lefschetz fibration can not be arbitrary and depends on the 
genus of the fiber, the genus of the base and the algebraic 
structure of the mapping class group of the generic fiber --- as we will
illustrate below.
  
Let $N(g,h)$ denote the minimal number of singular fibers in a
relatively minimal genus $g$ Lefschetz fibration (with at least
one singular fiber) over a closed (Riemann) surface of genus $h$.

Our first result is the following theorem.

\begin{thm} \label{thm1}
$N(g,h)=1$ iff $g \geq 3$ and $h \geq 2$.
\end{thm}

Let $\S _g$ be a closed orientable surface of genus $g$. The
mapping class group ${\CM}_{g}$ of $\S _g$ is defined to be the
group of isotopy classes of orientation preserving diffeomorphisms
$\S _g\to \S _g$. 

For a simple closed curve $a$ on an oriented surface, let us
denote by $t_a$ the right Dehn twist about $a$.

We also prove the following theorem, which is needed for the proof
of Theorem \ref{thm1}.

\begin{thm} \label{thm2}
Let $\S _g$ be a closed connected oriented surface of genus $g$
and let $a$ be a simple closed curve on $\S _g$. If $g\geq 3$,
then $t_a$ can be written as a product of two commutators.
\end{thm}

\section{Definitions and Proofs}

\begin{df} Let $M$ be a closed, connected, oriented smooth four manifold.
A Lefschetz fibration is a map $\pi : $M$ \rightarrow \S$, where $\S$ is a
closed, connected, oriented surface, such that $\pi$ is injective on the 
set of critical points $C= \{x_1,...,x_n \}$ and about each $x_i$ and 
$\pi (x_i)$ there are complex local coordinate charts agreeing with
the orientations of $M$ and $\S$ on which $\pi$ is of the form 
$$\pi (z_1 , z_2) = z_1^2 + z_2^2 .$$ 
\end{df}

Any fiber containing a critical point is called a singular fiber. 
Clearly all regular fibers are closed surfaces and are of the same 
diffeomorphism type. We will assume that the generic fiber is 
connected and our fibration is relatively minimal, i.e., no fiber 
contains a (-1) sphere. If the genus of the fiber is at least two,
then  a Lefschetz fibration is determined by its monodromy representation
$$ {\pi}_1 (\S - \{x_1,...,x_n \}) \rightarrow {\CM}_g . $$ 
(See \cite{gs} for more about Lefschetz fibrations.) 

Recall that for a group $G$, the first homology group with
integral coefficient is $H_1(G)=G/[G,G]$, where $[G,G]$ is the
commutator subgroup of $G$, the subgroup generated by the
commutators $[x,y]=xyx^{-1}y^{-1}$ for all $x,y\in G$. Let $\S_g$
be a closed oriented surface. It is well known that $H_1({\CM}_g)$
is generated by the class of a Dehn twist about a nonseparating
simple closed curve and is equal to ${\IZ}_{12}$ if $g=1$,
${\IZ}_{10}$ if $g=2$ and trivial if $g\geq 3$. We say that a
simple closed curve on a closed surface is nontrivial if it does
not bound a disc.

\begin{lemma} \label{lemma3}
 For every $h \geq 0$, $N(1,h)=12$.
\end{lemma}

{\it Proof.} All nontrivial circles on $\S_1$ are nonseparating
and the right Dehn twists map to the same generator under the
natural map ${\CM}_1\rightarrow H_1({\CM}_1)$. Hence, if a product
of right Dehn twists is trivial (or equal to a product of commutators), 
then the number of twists is
divisible by $12$. This implies that the number of singular fibers
in a genus one Lefschetz fibration must be divisible by $12$.
Moreover, this number is clearly realized by the elliptic surfaces
$ E(1) \# \S_1 \times \S_h$ for any $h \geq 0$.

\begin{prop} \label{prop4}
For every $h \geq 0$, $5\leq N(2,h) \leq 8$.
\end{prop}

{\it Proof.} Since every right Dehn twist about a nontrivial
separating simple closed curve on $\S _2$ is the product of $12$
right Dehn twists about nonseparating simple closed curves, the
order of its image under the map ${\CM}_2\rightarrow H_1({\CM}_2)$
is $5$. Hence, if a product of right Dehn twists about $n$
nonseparating and $s$ separating simple closed curves is trivial 
(or equal to a product of commutators), 
then $$n+2s \equiv 0 \;\; (\bmod 10).$$ Therefore, there are at
least $5$ singular fibers 
in a genus two Lefschetz fibration over a closed surface.
Moreover, a genus $2$ Lefschetz fibration with $8$ singular fiber
is constructed in \cite{m2}.

\bigskip

For a $4$-manifold $Y$, let $e(Y)$ and $\sigma (Y)$ denote the
Euler characteristic and the signature of $Y$, respectively. The
next proposition is proved by Stipsicz \cite{s1}.

\begin{prop}
Let $M$ be a 4-manifold which admits a relatively minimal genus $g
>0 $ Lefschetz fibration over a surface of nonzero genus. Then $0
\leq c_1^2 (M) \leq 10 \chi (M)$, where  $c_1^2 (M) =
3\sigma(M)+2e(M)$ and $\chi (M)= (\sigma(M) +e(M))/4 $.
\label{prop6}
\end{prop}

\begin{prop} \label{prop7}
$N(g,1) > 1$ for all $g \geq 1$.
\end{prop}
{\it Proof.}
 Suppose that $N(g,1) = 1$ for some $g$.  Then there is a
$4$-manifold $Y$ which admits a relatively minimal genus $g$
Lefschetz fibration over a torus with only one singular fiber.
Hence $e(Y)=1$. The inequality $ 0 \leq c_1^2 (Y) =
3\sigma(Y)+2e(Y)$ implies that $\sigma (Y) \geq -2/3$. This gives
$\sigma (Y) \geq 0$ since $\sigma (Y)$ is an integer.

On the other hand, $\chi (Y)= \frac{\sigma(Y) +e(Y)}{4}=
\frac{\sigma(Y) + 1}{4}$. Therefore $\sigma(Y)=4 \chi (Y) -1$. So
we get $c_1^2 (Y)= 3\sigma(Y)+2e(Y)=12 \chi (Y) -1 \leq 10 \chi
(Y)$. Hence $ \chi (Y) \leq 0$ and $\sigma(Y) \leq -1$ which is a
contradiction. This proves the proposition.

\begin{cor} \label{cor1}
The Dehn twist about a simple closed curve on a closed surface
cannot be equal to a single commutator.
\end{cor}
{\it Proof.}
 Suppose that a right Dehn twist is a commutator. Then there is a
relatively minimal Lefschetz fibration of genus $g$ with only one
singular fiber over the torus. Hence, $N(g,1)=1$, which
contradicts to Proposition \ref{prop7}.

\bigskip
We are now ready to prove our main results.

\bigskip
{\it Proof of Theorem \ref{thm2}.} Consider a sphere $X$ with four
holes with boundary components $a,a_1,a_2,a_3$. By the lantern
relation \cite{j} there are three simple closed curves
$b_1,b_2,b_3$ on $X$ such that
\[ t_at_{a_1}t_{a_2}t_{a_3}=t_{b_1}t_{b_2}t_{b_3} \]
or
\[ t_a=
t_{b_1}t_{a_1}^{-1}t_{b_2}t_{a_2}^{-1}t_{b_3}t_{a_3}^{-1}.\]
 Here for two diffeomorphisms $\a$ and $\b$, the composition $\a\b$ means
that $\a$ is applied first.

\begin{figure}[hbt]
\centerline{\psfig{figure=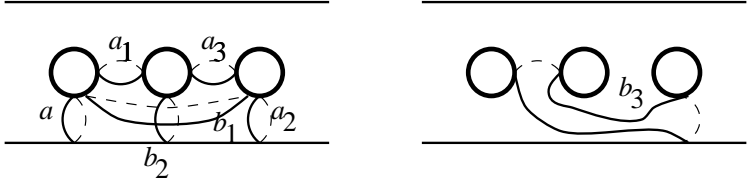,height=2.5cm,width=10cm}}
\caption{$a$ is nonseparating.} \label{sekil1}
\end{figure}

\begin{figure}[hbt]
\centerline{\psfig{figure=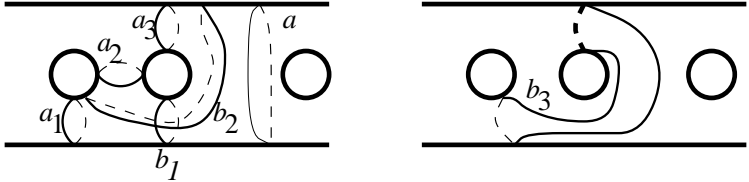,height=2.5cm,width=10cm}}
\caption{$a$ is separating.} \label{sekil2}
\end{figure}

Since the genus of $\S$ is at least three, $X$ can be embedded in
$\S$ in such a way that $a_1,a_2,a_3,b_1,b_2,b_3$ are all
nonseparating. The simple closed curve $a$ can be chosen either
nonseparating or separating bounding a subsurface of arbitrary
genus (cf. Fig.\ref{sekil1} and Fig.\ref{sekil2}). Furthermore,
the complement of $a_1\cup b_1$ and that of $a_2\cup b_2$ are
connected.  Hence, there is an orientation preserving
diffeomorphism $F$ of $\S$ such that $F(a_1)=b_2$ and
$F(b_1)=a_2$. Let $K$ be another orientation preserving
diffeomorphism of $\S$ such that $K(b_3)=a_3$ and let $f$ and $k$
be the isotopy classes of $F$ and $K$ respectively. Then

\begin{tabular}{ccl}
  $t_a$ &=& $t_{b_1}t_{a_1}^{-1}t_{F(a_1)}t_{F(b_1)}^{-1}
                                   t_{b_3}t_{K(b_3)}^{-1}$
                                    \vspace*{0.3cm}
                                   \\

        &=& $t_{b_1}t_{a_1}^{-1} f^{-1} t_{a_1}t_{b_1}^{-1}f
                            t_{b_3}k^{-1}t_{b_3}^{-1}k$
                            \vspace*{0.3cm} \\
        &=& $[t_{b_1}t_{a_1}^{-1}, f^{-1} ] [t_{b_3},k^{-1}].$
\end{tabular}\\

This concludes the proof of Theorem \ref{thm2}.

\bigskip

{\it Proof of Theorem $1$.} The fact that $N(g,h)>1$ is proved
above for $g < 3$ or $h = 1$. Clearly, $N(g,0)>1$ for all $g\geq
1$.

Suppose that $g \geq 3$. Using the relation $t_a =
[t_{b_1}t_{a_1}^{-1}, f^{-1} ] [t_{b_3},k^{-1}]$, we can construct
a genus $g$ Lefschetz fibration with only one singular fiber
(whose vanishing cycle is the given curve $a$) over a closed
surface of genus $2$ as follows:

The fundamental group of $\S _2$ has a presentation
 $$
    {\pi}_1(\S_2)= \langle x_1, y_1, x_2, y_2\, | \,[x_1, y_1][x_2,y_2]
    =1 \rangle.
 $$
 By the standard bundle theory, there exist a genus $g$
surface bundle over $\S_2$ with the following monodromy
representation:
 $$
    \Psi :{\pi}_1(\S_2) \rightarrow  {\CM}_{g},
 $$
 where $ \Psi (x_1) = t_{b_1}t_{a_1}^{-1} $, $ \Psi (y_1) = f^{-1} $,
$ \Psi (x_2) = t_{b_3} $, $ \Psi (y_2) = k^{-1}.$

Consequently there exists a genus $g$ surface bundle over the
surface ${\S_{2,1}}$ of genus two with one hole, such that the
monodromy over the boundary is given by the product
$$[t_{b_1}t_{a_1}^{-1}, f^{-1}] [t_{b_3},k^{-1}].$$

We can extend this surface bundle to obtain a Lefschetz
fibration over the closed genus two surface $\S_2$ just by
inserting a singular fiber with monodromy $t_a$. To obtain a
Lefschetz fibration with one singular fiber over a surface of
genus $h \geq 3$, we can take a fiber sum with $\S_g \times
\S_{h-2}$. This proves that $N(g,h)=1$ if $g \geq 3$ and $h \geq 2$.

This finishes the proof the theorem.

 \bigskip
 \noindent
 {\bf Remarks.}

1) Theorem \ref{thm2} answers a question of Mess in negative. See
Problem $2.13$ in \cite{k}. If the genus of the surface is at
least three, then the proof of the Theorem \ref{thm2} shows that
$t^n_a$ can be expressed as a product of $(3n+1)/2$ (resp. $3n/2$)
commutators if $n$ is odd (resp. even). Using the lantern
relation, we  can also prove that if the genus is at least four,
then $t_a^n$ can be written as a product of $n$ (resp. $n+1$)
commutators if $n$ is even (resp. odd).

2) Theorem \ref{thm2} also shows that the unique singular fiber
can be chosen reducible or irreducible. (A singular fiber is 
called reducible if the corresponding vanishing cycle is separating.)

3)ESTIMATES ON $N(g,0)$.

In \cite{c}, Cadavid proves that $N(g,0) \leq 2g+4$ if $g$ is even
and $N(g,0) \leq 2g+10$ if $g$ is odd. (This result was also
discovered independently by the authors \cite{ko}). Recently,
Stipsicz \cite{s2} proved that $\frac{1}{5} (4g+2) \leq N(g,0)$.

4)ESTIMATES ON $N(2,h)$.

In \cite{m2}, it is proven that  $5 \leq N(2,h) \leq 8$ for all $h
\geq 0$. The second author showed that $N(2,0)=7$ or $8$ in his
thesis \cite{o}. It is also clear that there is a $k \geq 1$ such
that $N(2,h)=5$ for each $h \geq k$ (cf. \cite{e}).

5) Examples of a Lefschetz fibration with a unique singular fiber
and prescribed fundamental group were constructed in \cite{abkp}.

\end{document}